\documentclass[12pt,leqno]{amsart}
\usepackage{amsmath}
\usepackage{graphicx,color}
\usepackage{verbatim}
\newtheorem{theorem}{Theorem}[section]

\newtheorem{lemma}[theorem]{Lemma}
\newtheorem{proposition}[theorem]{Proposition}
\theoremstyle{definition}
\newtheorem{definition}[theorem]{Definition}
\theoremstyle{remark}

\numberwithin{equation}{section}

\def\R{\mathbb{R}} 
\def\C{\mathbb{C}} 

\def\QH{\mathbb{H}} 

\def\rt{\rightarrow}

\begin{document}

\sloppy

\title[Gromov hyperbolicity of convex domains]{On the Gromov hyperbolicity of \\  convex domains in $\mathbb C^n$}

\author{Herv{\'e} Gaussier and Harish Seshadri}

\begin{abstract}
We give a necessary complex geometric condition for a bounded smooth convex domain in $\mathbb C^n$, endowed with the Kobayashi distance, to be Gromov hyperbolic. More precisely, we prove that if a $C^\infty$ smooth bounded convex domain in $\C^n$ contains an analytic disk in its boundary, then the domain is not Gromov hyperbolic for the Kobayashi distance. 

We also provide examples of bounded smooth convex domains that are not strongly pseudoconvex but are Gromov hyperbolic.
\end{abstract}

\maketitle

\section{Introduction and Main Result}

The notion of Gromov hyperbolicity (or `` $\delta$-hyperbolicty") of a metric space, introduced by M.Gromov in \cite{gro}, can be loosely described as ``negative curvature at large scales".   The prototype of a Gromov hyperbolic space is a simply connected complete Riemannian manifold with sectional curvature bounded above by a negative constant. One of the reasons for studying Gromov hyperbolic spaces is that such a space inherits many of the features of the prototype, even though the underlying space may not be a smooth manifold and the distance function may not arise from a Riemannian metric.

There are extensive studies of interesting classes of Gromov hyperbolic spaces in the literature, the class of word-hyperbolic discrete groups perhaps being the most studied. In a different direction, domains in Euclidean spaces endowed with certain natural Finsler metrics have also been analyzed from this point of view. For instance, M.Bonk, J.Heinonen and P.Koskela \cite{bo-he-ko} studied the Gromov hyperbolicity of planar domains endowed with the quasihyperbolic metric and
P.H\"ast\"o, H.Lind\'en, A.Portilla, J.M.Rodriguez and E.Touris~\cite{ha-li-po-ro-to} obtained the Gromov hyperbolicity of infinitely connected Denjoy domains, equipped either with the hyperbolic metric or with the quasihyperbolic metric, from conditions on the
Euclidean size of the complement of the domain.
In higher dimensions Z.Balogh and S.Buckley \cite{ba-bu} give conditions equivalent to the Gromov hyperbolicity for domains contained in $\overline{\mathbb R^n}$ and endowed with the inner spherical metric (or with the inner Euclidean metric if the domain is bounded). Y.Benoist \cite{be1,be2} gave among other results a necessary and sufficient condition, called quasisymmetric convexity, for a bounded convex domain of $\mathbb R^n$ endowed with its Hilbert metric to be Gromov hyperbolic.

In this paper we study domains in $\C^n$ endowed with the Kobayashi metric. The Kobayashi pseudodistance was introduced by S.Kobayashi as a tool to study geometric and dynamical properties of complex manifolds. A systematic study of its main properties and applications can be found in \cite{kob1,kob2,kob3}. This pseudodistance describes in a very precise way whether a complex manifold contains arbitrary large complex discs. 
The Kobayashi metric on domains in complex manifolds has proved to be a powerful tool in different problems such as the biholomorphic equivalence problem, or extension phenomena for proper holomorphic maps.

The class of strongly convex domains admits a particularly rich theory,  by the work of L.Lempert \cite{lem1,lem2,lem3}. L.Lempert proved that such a marked domain $(D,p)$ admits a singular foliation by complex geodesics discs, that these complex geodesics are the only complex one-dimensional holomorphic retracts and that the associated Riemann map is solution of a homogeneous complex Monge-Amp\`ere equation. As an application he exhibited a complete set of invariants for such marked strongly convex domains. See  \cite{bl-du} for related results.

The behaviour of real geodesics in a strongly convex domain endowed with the Kobayashi distance is well understood since one can prove that every such real geodesic is contained in a complex geodesic.
The behaviour of real geodesics on a more general bounded domain endowed with the Kobayashi distance is related to the Gromov hyperbocility of the given domain. The first and essentially the only result in that field is due to Z.Balogh and M.Bonk who proved in \cite{bal-bon} that every strongly pseudoconvex bounded domain in $\mathbb C^n$ endowed with the Kobayashi distance is Gromov hyperbolic. They also proved that the polydisc, which is complete Kobayashi hyperbolic, is not Gromov hyperbolic.

The Balogh-Bonk result naturally raises the question of Gromov hyperbolicity  of weakly pseudoconvex domains in $\C^n$. As mentioned above (see also the remarks following Theorem \ref{nonhyp-thm}) it is easy to prove that the polydisc is not Gromov hyperbolic. However its boundary is not smooth. While it appears reasonable to expect that the same conclusion should hold for smooth weakly pseudoconvex bounded domains, such a result was not known, to the best of our knowledge. In this paper we prove that certain smooth bounded weakly pseudoconvex (but not strongly pseudoconvex) domains are not Gromov hyperbolic and, perhaps more surprisingly, some are.

More precisely, we investigate the Gromov hyperbolicity of smooth bounded convex domains endowed with the Kobayashi metric. We prove that the existence of a non constant analytic disc in the boundary of the domain is an obstruction to the Gromov hyperbolicity, giving a general answer to a question raised by S.Buckley.
We also provide examples of Gromov hyperbolic convex domains of finite type, in the sense of D'Angelo, that are not strongly pseudoconvex.

\vspace{2mm}
The precise statements of our results are as follows. Let $(D,d)$ be a metric space. A curve $\gamma: [a,b] \rt D$ is a {\it geodesic} if $\gamma$ is an isometry for the usual distance function on $[a,b] \subset \R$, i.e., $d(\gamma(t_1), \gamma_(t_2))= \vert t_1-t_2 \vert$ for all $t_1,t_2 \in [a,b]$. A {\it geodesic triangle} in $D$ is a union of images of three geodesics $\gamma_i: [a_i,b_i] \rt D$, $i=1,2,3$,  such that $\gamma_i(b_i)=\gamma_{i+1}(a_{i+1})$ where the indices are taken modulo $3$.  The image of each $\gamma_i$ is called a {\it side} of the geodesic triangle.  $(D,d)$ is
{\it Gromov hyperbolic} or $\delta$-hyperbolic if there exists $\delta >0$ such that for any geodesic triangle in
$D$ the image of every side is contained in the $\delta$-neighbourhood of the other two sides.

If $D \subset \C^n$ is a bounded domain, we denote the Kobayashi metric on $D$ by $d_D^K$.

The main result of this paper is the following

\begin{theorem}\label{nonhyp-thm}
Let $D$ be a bounded convex domain in $\mathbb C^n$ with a $\mathcal C^\infty$ boundary ${\partial D}$. Assume that ${\partial D}$ contains an analytic disk. Then $(D,d_D^K)$ is not Gromov hyperbolic.
\end{theorem}

The second result is an observation that certain weakly convex domains are Gromov hyperbolic:

\begin{theorem}\label{CompEllip-thm}
For $p \ge 1$ the {\it complex ellipsoid $D_p$ in} $\C^n$ given by
$$D_p = \{(z',z_n) \in \C^{n-1} \times \mathbb C \ \vert \ \vert z' \vert^2 + \vert z_n \vert^{2p} < \ 1\},$$
is Gromov hyperbolic for the Kobayashi distance.
\end{theorem}

In light of these results, an interesting problem consists in giving a precise complex geometric characterization of bounded smooth convex domains that are Gromov hyperbolic. In particular, we have the following question: \vspace{2mm}

{\it Is the Gromov hyperbolicity of a $C^\infty$ smooth bounded convex domain $D$ in $\C^n$ equivalent to the condition that $D$ is of finite type in the sense of D'Angelo?}\vspace{2mm}

We make a few remarks about the proofs of these results briefly, beginning with Theorem \ref{nonhyp-thm}.  It is an easy observation that the product of two complete noncompact geodesic metric spaces endowed with the maximum metric is not Gromov hyperbolic. We can see this for the Kobayashi metric on the bidisc $D=\Delta^2 \subset \mathbb C^2$ as follows. Consider the origin $O$ and the points $p_m:(1-\frac{1}{m},-1+\frac{1}{m}), q_m:(-1+\frac{1}{m},-1+\frac{1}{m})$. Denote by $s_m^1$ (resp. $s_m^2$) the unique geodesic, for the Kobayashi distance on $D$, joining $O$ to $p_m$ (resp. to $q_m$). Denote by $l_m$ the Kobayashi length of the unique geodesic joining $p_m$ to $q_m$. Then the unique point $z_m$ such that $d_{\Delta^2}^K(p_m,z_m) = \frac{l_m}{2}$ is at a Kobayashi distance $d_m$ from the two geodesics $s_m^1$ and $s_m^2$, with $\lim_{m \rightarrow \infty}d_m=+\infty$. Hence $(\Delta^2,d_{\Delta^2}^K)$ is not Gromov hyperbolic. Our proof of the non Gromov hyperbolicity of the domain $D$ in Theorem~\ref{nonhyp-thm} is inspired by that construction.

Of course, the main difficulty in dealing with a domain $D$ which is not a product is describing the geodesics (or quasi-geodesics). On the one hand, it is reasonable to expect that the metric behaves (in terms of geodesics) like a product near the holomorphic disc $H \subset \partial D$. On the other, the smoothness of $\partial D$ forces strict convexity, at some points close to $\partial H$, of $\partial D$. One of the main technical points of this work  is that the metric does behave like a product and the the aspect mentioned in the second point above does not dominate.

For our second result, we use results of K. Azukawa - M. Suzuki and J. Bland stating that  the sectional curvatures of the Bergman and K\"ahler-Einstein metrics on complex ellipsoids are bounded above by negative
constants. In particular, these metrics are Gromov hyperbolic. One can then use the Schwarz lemma in various forms to compare these metrics with the Kobayashi metric and conclude that it is Gromov hyperbolic as well.

\section{Notations and Definitions}
$\bullet$ For $1 \le i \le n$, $\pi_i: \C^n \rt \C$ will denote the $i$-th projection. \vspace{2mm}

$\bullet$ For $1 \le i \le n$, $Z_i= \{z \in \C^n \ \vert \ \pi_j(z)=0$ if $i \neq j \}$, $X_i = Z_i \cap \R^n$ and
$Y_i=Z_i \cap \sqrt{-1} \R^n$.

We will refer to these as the $z_i \ axis$, etc.   \vspace{2mm}

$\bullet$ If $z$ and $z'$ are two points in $\mathbb C^k$ we will denote by $dist_{eucl}(z,z')$ the Euclidean distance between $z$ and $z'$. \vspace{2mm}

Let $E$ be a smooth bounded convex domain in $\C^k$ and $z \in E$. \vspace{2mm}

$\bullet$ $\delta_E(z)$ will denote the Euclidean distance from $z$ to the boundary ${\partial E}$ of $E$. \vspace{2mm}

$\bullet$ For $v \in \mathbb C^k \backslash\{0\}$ we denote by $\delta_E(z,v)$ the Euclidean distance from $z$ to the boundary ${\partial E}$ of $E$ along the complex line $L(z,v):=\{z+\lambda v,\ \lambda \in \mathbb C\}$.\vspace{2mm}

$\bullet$ If $l \subset \C^n$ is a real line, $\delta_E(z,l)$ will denote the Euclidean distance from $z$ to ${\partial E}$ along $z+l$. \vspace{3mm}

$\bullet$ For $q \in \partial E$ let $r_q$ be the largest real number $r$ such that there is a ball $B$  of radius $r$ contained in $E$ with $\partial B$ tangent to $\partial E$ at $q$. \vspace{2mm}

We consider the complex Euclidean space $\mathbb C^n$ with coordinates $(z_1,\dots,z_n)$. For a positive real number $r$ we denote by $\Delta_{r}$ the disk $\Delta_r:=\{\zeta \in \mathbb C : |\zeta| < r\}$.

\begin{definition}\label{hyp-def}
Let $D$ be a domain in $\mathbb C^n$. \vspace{2mm}

$(i)$ The Kobayashi infinitesimal pseudometric is defined on $TD=D \times \mathbb C^n$ by :
$$
\forall (p,v) \in D \times \mathbb C^n,\ K_D(p,v):=\inf\{\alpha > 0/\ \exists \ f : \Delta \rightarrow D,\ f \ {\rm holomorphic}, \ f(0)=p,\ f'(0)=v/\alpha\}.
$$

$(ii)$ The Kobayashi peudodistance $d_D^K$ is defined on $D \times D$ by :
$$
\forall (p,q) \in D \times D,\ d_D^K(p,q) = \inf\left\{\int_0^1 K_D(\gamma(t),\gamma'(t))dt\right\}
$$
where the infimum is taken over all $\mathcal C^1$-paths from $[0,1]$ to $D$ satisfying $\gamma(0)=p,\ \gamma(1)=q$.

$(iii)$ The domain $D$ is {\sl Kobayashi hyperbolic} if $d_D^K$ is a distance on $D$.

$(iv)$ The domain $D$ is {\sl complete hyperbolic} if $(D,d_D^K)$ is a complete metric space.
\end{definition}

\begin{definition}\label{quasi-def}
Let $D$ be a domain in a manifold $M$ endowed with a continuous Finsler metric $K: TM \rightarrow \mathbb R$.  Let $\gamma : [a,b] \rightarrow E$ be a path joining two points $p$ and $q$ in $E$. \vspace{2mm}

(i) If $\gamma$ is piecewise smooth, the {\it length} of $\gamma$ is the quantity
$$l^K_D(\gamma) = \int_a^b K(\gamma(t), \gamma'(t))dt.$$ \vspace{2mm}

(ii) The {\it distance} between $p$ and $q$ is given by
$$d^K_D(p,q)= \inf \ l^K_D(\gamma)$$
where the infimum is taken over all piecewise smooth curves between $p$ and $q$. \vspace{2mm}

(i) $\gamma$ is a {\it geodesic} if $\gamma: [a,b] \rt D$ is an isometry for the usual distance function on $[a,b] \subset \R$, i.e., $d(\gamma(t_1), \gamma_(t_2))= \vert t_1-t_2 \vert$ for all $t_1,t_2 \in [a,b]$.  \vspace{2mm}

(ii) Let $A \ge 1$ and $B>0$. We say that $\gamma$ is a $(A,B)$ {\sl quasi-geodesic} if for every $t_1, t_2 \in [a,b]$ we have
$$
\frac{1}{A}\vert t_1- t_2 \vert -B \leq  d(\gamma(t_1),\gamma(t_2)) \leq A \vert t_1 -t_2 \vert + B.
$$ \vspace{2mm}

(iii) A {\it geodesic triangle} in $D$ is a union of images of three geodesics $\gamma_i: [a_i,b_i] \rt D$, $i=1,2,3$,  such that $\gamma_i(b_i)=\gamma_{i+1}(a_{i+1})$ where the indices are taken modulo $3$.  The image of each $\gamma_i$ is called a {\it side} of the geodesic triangle.  $(D,d_D^K)$ is
{\it Gromov hyperbolic} or $\delta$-hyperbolic if there exists $\delta >0$ such that for any geodesic triangle in
$D$ the image of every side is contained in the $\delta$-neighbourhood of the other two sides.
\end{definition}

\section{Proof of Theorem~\ref{nonhyp-thm}}
The proof of Theorem~\ref{nonhyp-thm} is a direct consequence of the following Proposition:
\begin{proposition}\label{flat-prop}
If $D$ satisfies the assumptions of Theorem~\ref{nonhyp-thm} then there exists $0 < r_0$ such that $D$ satisfies the following conditions: \vspace{2mm}

\noindent ${\bf (1)}$ $0 \in {\partial D}$, the tangent space to $D$ at $0$ is given by $T_0(D)=\{z \in \mathbb C^n : Re(z_n) = 0\}$ and  $D \subset \{z=(z',z_n) \in \mathbb C^{n-1} \times \mathbb C : Re(z_n) > 0\}$.  Moreover there is a
convex set $C \subset T_0(D) \cap \overline D$ containing $0$ such that $C$ is an open subset of $Z_1$.

\vspace{2mm}

 Assume that $C \cap X_1 = \{te_1\ \vert \ -a < t < a \}$ for some $a >0$. Let $\pm A= ( \pm a,0,...,0)$, $p^r= (-a,0,...,0,r)$,  $q^r=(a,0,...,0,r)$. \vspace{2mm}

\noindent ${\bf (2)}$ Let $K_r^\pm$ be the real line segments joining $q^r$ to $A$ and
$p^r$ to $-A$ respectively.  If $0< r<\sqrt{r_0}$ we have $$\delta_{D \cap (\pm A+Z_n)}(z) = Re(z_n)$$ for any $z \in K_r^\pm$ .
\vspace{2mm}

\noindent ${\bf (3)}$ There exists $\alpha >0$ with the following property: Let $E_r$ be the two-dimensional convex set $E_r= (Z_1+ (0,...,0,r)) \cap D$ and $L_r$  the real line segment joining $p^r$ to $q^r$. Then
$$  \delta_{E_r} (z) \ge \alpha  \ \delta_{E_r} (z, L_r) $$
for all $z \in L_r$ and $0<r <\sqrt{r_0}$.

\vspace{2mm}

\noindent ${\bf (4)}$ If $0 < r_1 < r_2 < \sqrt{r_0}$ and $(z',r_1) \in D$ for some $z' \in \C^{n-1}$ then $(z',r_2) \in D$.
\vspace{2mm}

\noindent ${\bf (5)}$ Assume $0<r<\sqrt{r_0}$. There exist $\beta>0$ and   smooth functions $f, \ h :[0,r_0] \rt [0,\infty)$  such that \vspace{2mm}

$(v_a)$ $\beta^{-1} h(r) \le \delta_{E_r}(p^r) \le \beta h(r)$ \ \ \ \ $\beta^{-1} f(r) \le \delta_{E_r}(q^r) \le \beta f(r)$

\vspace{2mm}

$(v_b)$ The functions $h$ and $f$ satisfy the following Condition, called Condition(**): \vspace{2mm}

\vskip 0,1cm
\noindent Condition(**) {\sl ``$f=g^{-1}$ (resp. $h=k^{-1}$) where $g$ (resp. $k$) is a strictly increasing  convex function of class $\mathcal C^\infty$ defined on $[0,\varepsilon]$ for some $\varepsilon > 0$ and such that $g^{(l)}(0) = 0$ (resp. $k^{(l)}(0)=0$) for every nonnegative integer $l$''.}

\end{proposition}

\noindent{\bf Proof of Proposition~\ref{flat-prop}:} Suppose that $p \in {\partial D}$ lies in the analytic disc $S$ in the hypothesis of Theorem ~\ref{nonhyp-thm}.  \vspace{2mm}

{\bf Claim:} $S \subset T_pD \cap D$.

Choose a complex affine linear map $T: \C^n \rt \C^n$ such that $T(p)=(0,...,0,1)$ and $T(T_pD)= \{ z \in \C^n \ \vert \ Re(z_n)=1 \}$. By considering
$-T$ if necessary, we can assume that $T(D) \subset \{ z \in \C^n \ \vert \ Re(z_n) >1 \}$, since $D$ is convex.
 Let $\phi: \Delta \rt S \subset {\partial D}$ be a holomorphic map with image $S$ with $T \circ \phi=((T \circ \phi)_1,...,(T \circ \phi)_n)$. If the holomorphic map $(T \circ \phi)_n: \Delta \rt \C$  is not constant then it violates the maximum modulus principle since
 $\vert (T \circ \phi)_n (z) \vert \ge 1$ for all $z \in \Delta$ and $\vert (T \circ \phi)_n(0)\vert =1$. Since $(T \circ  \phi)_n(0)=1$ we get
 $(T \circ \phi)_n(z)=1$ for all $z \in \Delta$. This proves the claim.\vspace{2mm}

Let $D_1 = T(D)- (0,...,0,1)$ where $T$ is as above. Note that $C=T_0D_1 \cap D_1$ is a convex subset of $\C^n$ containing $S$. We now
invoke the easy fact that {\it given a convex subset $C$ of Euclidean space containing the origin there is a linear subspace $V$ containing $C$ as an open (in $V$) subset}. Take any point $0$ in $S_1$, where $S_1$ is the holomorphic disc in $D_1$ corresponding to $S$ in $D$. By the above fact, if  $L$ is the tangent space to $S_1$ at $0$, then $E=L \cap C$ is open in $L$ and convex. Since $L$ is a complex line we can find  a complex rotation $R: \C^n \rt \C^n$
with  $R(L)=Z_1$ and $R(e_{2n-1})=e_{2n-1}$. We still have $0 \in R(D_1)$, $R(D_1)\subset \{ z \in \C^n \ \vert \ Re(z_n) >0 \}$ and $T_0R(D_1)=\{ z \in \C^n \ \vert \ Re(z_n) =0 \}$. Since $C$ is convex so is $R(C)$.   This ensures Condition ({\bf 1}) is satisfied.
\vskip 2mm
We relabel $R(D_1)$ as $D$. We will need the following simple lemma, the proof of which we will skip, for the other conditions:

\begin{lemma}\label{kml}
Let $E$ be a bounded convex domain with ${\mathcal C}^2$ boundary in $\C$ and $L$  an affine real line in $\C$.

If $p \in {\partial E} \cap L$ and this intersection is transversal at $p$, then there is a neighbourhood $U$ of $p$ in $E$  such that
 $$  \delta_{E} (z) \ge  cos(\theta) \ \delta_{E} (z, L) $$
for all $z \in L \cap U$, where $\theta$ is the angle between $L$ and the inward normal to ${\partial D}$ at $p$. If $L$ is perpendicular to ${\partial E}$ at $p$, we can assume that $\delta_{E} (z) =  \delta_{E} (z, L)$ for all $z \in L \cap U$.

We can take $U=B(p,\epsilon) \cap E$ where $\epsilon$ depends only on $r_p$ and $\theta$.
\end{lemma}

We come to Condition ({\bf 2}): Let $S$ be the two-dimensional convex set $S= D \cap (Z_n+A)$. By the convexity of $S$ and the fact that the tangent space to $\partial S$ at $(a,0,...,0)$ is the axis $Im(z_n)$, Lemma \ref{kml} gives a neighbourhood $U$ of $(a,0,...,0)$ in $S$ on which
$$\delta_{D \cap (\pm A+Z_n)}(z) = Re(z_n)$$ for any $z \in U \cap K^+_r$. A similar statement is true for $p^r$. Hence there exists $r_1 >0$ such that Condition ({\bf 2}) is satisfied for $r <r_1$.  \vspace{2mm}

To see the validity of Condition ({\bf 3}), the last part of Lemma \ref{kml} implies that there exists $0<r_2 < r_1$, $\alpha >0$ and $\gamma >0$
such that
$$  \delta_{E_r} (z) \ge \alpha  \ \delta_{E_r} (z, L_r) $$
for all $z \in L^{\pm \delta}_r$ and all $0<r<r_2$, where $L^{\pm \delta}_r$ are the real segments connecting $(-a+\gamma,0,...,0,r)$
to $p^r$ and $q^r$ to $(a -\gamma,0,...,0,r)$ respectively.  For the remaining parts of the $L_r$ note that the function $\phi: [-A+\gamma,\ A+\gamma] \times [0,\ r_2] \rt [0,\infty)$ given by
$\phi(s,r) = \frac {\delta_{E_r} (z(s,r))}{\delta_{E_r} (z(s,r), L_r)}$ is continuous, where
$z(s,r) = (s,0,...,0,r)$.  Hence its infimum is attained and, in particular, positive. \vspace{2mm}

Condition ({\bf 4}) easily follows from the compactness of $\overline E$.  \vspace{2mm}

For $0<r<r_2$ let $f(r)=\delta_{E_r} (q^r, L_r)$ and $h(r)=\delta_{E_r} (p^r, L_r)$. Condition $(v_a)$ now follows from Condition ({\bf 3}). Note that there exists $0<r_3<r_2$ such that the plane curve ${\partial D} \cap \{z \in \C^n \ \vert \ z=(s,0,...,0,r), \ a \le s, \ 0\le r \le r_3 \}$ (resp. ${\partial D} \cap \{z \in \C^n \ \vert \ z=(s,0,...,0,r), \ s \le -a, \ 0\le r \le r_3 \}$) is a graph of  a strictly increasing (resp. strictly decreasing)convex $C^\infty$ function. These functions are precisely
$x \rt g(x-a)$ and $k(-x-a)$.  \vspace{2mm}

Finally we take $r_0= r_3$.

\qed.

We can now prove Theorem~\ref{nonhyp-thm}. Assume by contradiction that $D$ is $\delta$-hyperbolic for some $\delta>0$. First we note that $(D, d_K^D)$ is a {\it geodesic} metric space, i.e. any two points in $D$ are connected by a geodesic. For a strongly convex domain this follows from  a basic result of L. Lempert \cite{lem1} which asserts that, in fact, there is a complex geodesic containing any two given points.
A complex geodesic in $D$ is an isometric map of $(\Delta, d_\Delta^K)$, the unit disc in $\C$ with the Poincare distance, into $(D,d^K_D)$. For a weakly convex domain the existence of complex geodesics containing any two points is due to H.L.Royden and B.Wong \cite{ro-wo} and a proof can be found in \cite{ab}, Theorem 2.6.19 p.265.

Let $A$ and $B$ be two positive constants. It follows from \cite{ghys-harpe} that there exists a constant $M>0$, depending only on the $\delta$ and on the constants $A,B$, such that the Hausdorff distance between any $(A,B)$ quasi-geodesic in $D$ and any geodesic between the endpoints of the quasi-geodesic is less than $M$. Hence in order to contradict our hypothesis of $\delta$-hyperbolicity it is sufficient to find two positive constants $A$ and $B$  and a family of $(A,B)$ quasi-triangles, namely unions of $(A,B)$ quasi-geodesics  joining three points, which violate the $\delta$-hyperbolicity condition for any $\delta$. More precisely we prove in Proposition~\ref{nonhyp-prop} that there exists $r_0>0$, positive constants $A, B$ and a point $x^0 \in D$ such that for every $0 < r < r_0$ there are points $p^r$ and $q^r$ in $D$ for which the curves $l^{x^0,p^r}$, $l^{x^0,q^r}$ and $\gamma^{r,r'}$ are $(A,B)$ quasi-geodesics. Finally we prove that for every $0 < r < r_0$ there is a point $z^{r,r'}$ on $\gamma^{r,r'}$ such that:
$$
\lim_{r \rightarrow 0}d^K_D(z^{r,r'},l^{x^0,p^r} \cup l^{x^0,q^r}) = +\infty.
$$
That condition violates the $\delta$-hyperbolicity condition as claimed.

\qed

\section{The Kobayashi metric on convex domains}
\subsection{Estimates for the Kobayashi metric}
We start with some general facts concerning convex domains in $\mathbb C^n$.
The following estimate of the Kobayashi infinitesimal pseudometric, obtained by I.Graham \cite{gra} and S.Frankel \cite{fra} will be an essential tool in the paper. See \cite{bed-pin} for an elementary proof.

\vskip 0,2cm
\noindent{\bf Proposition A.}
{\sl Let $D \subset \mathbb C^n$ be a convex domain. If $a \in D$ and if $v$ is a tangent vector, then
$$
\frac{|v|}{2\delta_D(a,v)} \leq K_D(a,v) \leq \frac{|v|}{\delta_D(a,v)}
$$
where $\delta_D(a,v)$ denotes the Euclidean distance from $a$ to $L(a,v) \cap {\partial D}$. Here $L(a,v)$ denotes the complex line passing through $a$ in the direction $v$.
}

\vskip 0,2cm
We will also use the following {\sl Boxing Lemma}:
\begin{lemma}\label{boxing-lem}
Let $R_{\alpha,\beta}:=\{\zeta \in \mathbb C : 0 < Re(\zeta) < \alpha,\ -\beta < Im(\zeta) < \beta\}$ and let $D_{R,\alpha,\beta}$ be the domain in $\mathbb C^n$ defined by $D_{R,\alpha,\beta}:=\Delta_R^{n-1} \times R_{\alpha,\beta}$, where $\alpha,\beta,R>0$.
Then for every $0 < r < r' < \alpha$ and for every $z^r=(z^r_1,\dots,z^r_n), \ z^{r'}=((z^{r'}_1,\dots,z^{r'}_n) \in D_{R,\alpha,\beta}$ with $Re(z^r_1)=r$, $Re(z^{r'}_1)=r'$ we have:
$$
d^K_{D_{R,\alpha,\beta}}(z^r,z^{r'}) \geq \frac{1}{2} \left|\log\left(\frac{r'}{r}\right)\right|.
$$
\end{lemma}

\vskip 0,2cm
\noindent{\bf Proof of Lemma~\ref{boxing-lem}.} Denote by $\pi_n$ the holomorphic projection from $D_{R,\alpha,\beta}$ to $R_{\alpha,\beta}$.
Let $\gamma=(\gamma_1,\dots,\gamma_n):[t_0,t_1] \rightarrow D_{R,\alpha,\beta}$ be a curve with $\gamma(t_0)=z^r$ and $\gamma(t_1)=z^{r'}$.
Since $\gamma_n=\pi_n(\gamma)$ we have from the Schwarz Lemma and from Proposition A:
$$
\begin{array}{lllll}
l^K_{D_{R,\alpha,\beta}}(\gamma) & \geq & l^K_{R_{\alpha,\beta}}(\gamma_n) & \geq & \displaystyle \frac{1}{2}\int_{t_0}^{t_1}\frac{|\gamma_n'(t)|}{\delta_{R_{\alpha,\beta}}(\gamma_n(t))}dt\\
& & & \geq & \displaystyle \frac{1}{2}\int_{t_0}^{t_1}\frac{Re(\gamma_n'(t))}{|Re(\gamma_n(t))|}dt
\end{array}
$$
since for every $t \in [t_0,t_1]$ we have $\delta_{R_{\alpha,\beta}}(\gamma_n(t)) \leq |Re(\gamma_n(t))|$ and $|\gamma_n'(t)| \geq Re(\gamma_n(t)) > 0$.
Finally:
$$
\int_{t_0}^{t_1}\frac{Re(\gamma_n'(t))}{|Re(\gamma_n(t))|}dt \geq
\left|\int_{t_0}^{t_1}\frac{Re(\gamma_n'(t))}{Re(\gamma_n(t))}dt\right| =
\left|\log\left(\frac{r'}{r}\right)\right|.
$$
\qed

\subsection{Quasi-geodesics in convex domains}
The following proposition provides very simple examples of quasi-geodesics, for the Kobayashi metric, in any bounded convex domain of class $\mathcal C^1$  in $\mathbb C^n$.
\begin{proposition}\label{quasi1-prop}
Let $D$ be a bounded convex domain in $\mathbb C^n$ of class $\mathcal C^1$ and let $x \in D$. Then every real segment $l$, parametrized with respect to Kobayashi arc-length, joining $x$ to $q \in {\partial D}$ is a $(A,B)$ quasi-geodesic where $A, \ B$ depend only on $x$, the angle between $l$ and the inward normal to $D$ at $q$ and $D$.
\end{proposition}

\noindent{\bf Proof of Proposition~\ref{quasi1-prop}.} For $q \in \partial D$ let $r_q$ be as in Section 2. By the continuity of the map $q \mapsto r_q$ and the compactness of $\partial D$, if $\alpha = \inf_{q \in \partial D} r_q$ then $\alpha >0$. Let $U = B_q(\frac {\alpha}{2})$.

 Let $l(t_0) \in \partial B_q(\frac {\alpha}{2})$ and $B_1 = l^K_D(l \vert_{[0,t_0]})= t_0$. Note that $B_1$ is bounded above by a constant depending only on $x$ and $D$. It is enough to prove that $l \vert _{[t_0, \infty)}$ is a quasi-geodesic as in Proposition \ref{quasi1-prop} with constants $(A, B_2)$. To see this let $0<t_1 <t_0< t_2  $. Then
 $$  d^K_D(l(t_1), l(t_2)) \le  d^K_D(l(t_0), l(t_2)) +B_1  \le A \vert t_0 -t_2 \vert +B_1+B_2 \le A \vert t_1 -t_2 \vert +B $$
where $B=B_1+B_2$.
Also
$$d^K_D(l(t_1), l(t_2)) \ge d^K_D(l(t_0), l(t_2)) \ge A^{-1} \vert t_0 -t_2 \vert -B_2 \ge A^{-1} \vert t_1 -t_2 \vert -B.$$
The other cases ($t_1 <t_2 <t_0$ and $t_0 < t_1<t_2$) can be seen similarly.

The rest of the proof is devoted to showing that $l \vert_{[t_0,\infty)}$ is a quasi-geodesic. By modifying $D$ by a complex affine linear map if necessary, we may assume that $q=0$ and that the tangent space $T_q({\partial D})$ is given by $T_q({\partial D})=\{(z',z_n) \in \mathbb C^{n-1}\times \mathbb C : Re(z_n)=0\}$. Denote the inward pointing normal vector to ${\partial D}$ at $q$ by $\nu$ and let $l_\nu = \R^+ \nu$. Note that $\nu=e_{2n-1}$.

  Let $p=(p',p_n) \in l \cap U$. There exists $p_\nu \in l_\nu$ with $Re({(p_\nu)}_n)=Re(p_n)$. Let $\gamma : [0,1] \rightarrow U \subset \mathbb C^n$ be the straight line joining $p$ to $p_\nu$. 
We note the following elementary fact, the proof of which we skip: There exists $C>0$ depending only on
$\alpha$ and $\theta$ such that $$\delta_D (\gamma (t), \gamma'(t)) \ge C  \delta_D(p_\nu)$$.
Hence we have
$$l^K_D(\gamma) \leq \int_0^1\frac{|\gamma'(t)|}{\delta_D(\gamma(t),\gamma'(t))}dt
\leq \frac{1}{C}\int_0^1 \frac{|\gamma'(t)|}{\delta_D(p_\nu)}dt =  \frac{1}{C} \frac{|p-p_\nu|}{\delta_D(p_\nu)}.
$$

Now, since $l$ makes an angle $\theta$ with $\nu$ we have $|p-p_\nu| =  \tan (\theta) \delta_D(p_\nu)$

Finally:
\begin{equation}\label{bounded-eq}
l_D^K(\gamma) \leq \frac{\tan (\theta)}{C}.
\end{equation}






In light of Equation~(\ref{bounded-eq})  it is sufficient to prove the following lemma to complete
the proof of Proposition~\ref{quasi1-prop}.

\begin{lemma}\label{quasi1-lem}
The real line $l_\nu \cap U$ is a $(1, log2)$ quasigeodesic for the Kobayashi distance on $D$.
\end{lemma}

\noindent{\bf Proof of Lemma~\ref{quasi1-lem}.} We continue with our assumption that $q=0$.
There exists $R>0$ such that $D \subset P=\Delta_R^{n-1} \times \QH$ where $\QH:=\{\zeta \in \mathbb C :  Re(\zeta) >0\}$.

Let $p^1,\ p^2$ be two points in $l_\nu \cap U$.
As before, write $(z',z_n) \in \mathbb C^{n-1} \times \mathbb C$ the coordinates corresponding to the decomposition of the product $P$, and denote $p^1=((p^1)',p^1_n)$, $p^2=((p^2)',p^2_n)$. We suppose that $|p^1_n| > |p^2_n|$.
We have :
$$
d^K_D(p^1,p^2) \geq d^K_P(p^1,p^2) \geq d^K_{\{(p^1)'\} \times \QH}(p^1,p^2).
$$
The first inequality comes from the decreasing property of the Kobayashi distance and the second inequality comes from the fact that the Kobayashi distance on a product domain is the maximum of the Kobayashi distances on each factor.

On the other hand
$$
d^K_{\{(p^1)'\} \times \QH}(p^1,p^2) \ge  d^K_{\{(p^1)'\} \times \Delta(p^1_n,|p^1_n|)}(p^1,p^2) -log 2.
$$
However since $\{(p^1)'\} \times \Delta(p^1_n,|p^1_n|) \subset D$ we have by the decreasing of the Kobayashi distance~:
$$
d^K_{\{(p^1)'\} \times \Delta(p^1_n,|p^1_n|)}(p^1,p^2) \geq d^K_D(p^1,p^2).
$$
If $\gamma:[0,1] \rt \QH$ is the straight line segment joining $p^1$ and $p^2$, then
$$
l^K_{\{(p^1)'\} \times \Delta(p^1_n,|p^1_n|)}(\gamma_n) = d^K_{\{(p^1)'\} \times \Delta(p^1_n,|p^1_n|)}(p^1,p^2).
$$
Hence
$$
 d^K_D(p^1,p^2) + log2 \geq l^K_D(\gamma_n) \geq d^K_D(p^1,p^2)
$$
this proves Lemma~\ref{quasi1-lem} and Proposition~\ref{quasi1-prop}. \qed

\section{Quasi-geodesics in ``flat'' convex domains}\label{flat-section}
In this section we focus on some special convex domains in $\mathbb C^n$ for which we construct specific quasi-geodesics for the Kobayashi metric. We assume that there exists $0 < r_0 < 1$ such that the domain $D$ satisfies the following five conditions:

\vskip 0,3cm

\noindent ${\bf (1)}$ $0 \in {\partial D}$, the tangent space to $D$ at $0$ is given by $T_0(D)=\{z \in \mathbb C^n : Re(z_n) = 0\}$ and  $D \subset \{z=(z',z_n) \in \mathbb C^{n-1} \times \mathbb C : Re(z_n) > 0\}$.  Moreover there is a
convex set $C \subset T_0(D) \cap \overline D$ containing $0$ such that $C$ is an open subset of $Z_1$.

\vspace{2mm}

 Assume that $C \cap X_1 = \{te_1\ \vert \ -a < t < a \}$ for some $a >0$. Let $\pm A= ( \pm a,0,...,0)$, $p^r= (-a,0,...,0,r)$,  $q^r=(a,0,...,0,r)$. \vspace{2mm}

\noindent ${\bf (2)}$ Let $K_r^\pm$ be the real line segments joining $q^r$ to $A$ and
$p^r$ to $-A$ respectively.  If $0< r<\sqrt{r_0}$ we have $$\delta_{D \cap (\pm A+Z_n)}(z) = Re(z_n)$$ for any $z \in K_r^\pm$ .
\vspace{2mm}

\noindent ${\bf (3)}$ There exists $\alpha >0$ with the following property: Let $E_r$ be the two-dimensional convex set $E_r= (Z_1+ (0,...,0,r)) \cap D$ and $L_r$  the real line segment joining $p^r$ to $q^r$. Then
$$  \delta_{E_r} (z) \ge \alpha  \ \delta_{E_r} (z, L_r) $$
for all $z \in L_r$ and $0<r <\sqrt{r_0}$.

\vspace{2mm}

\noindent ${\bf (4)}$ If $0 < r_1 < r_2 < \sqrt{r_0}$ and $(z',r_1) \in D$ for some $z' \in \C^{n-1}$ then $(z',r_2) \in D$.
\vspace{2mm}

\noindent ${\bf (5)}$ Assume $0<r<\sqrt{r_0}$. There exist $\beta>0$ and   smooth functions $f, \ h :[0,r_0] \rt [0,\infty)$  such that \vspace{2mm}

$(v_a)$ $\beta^{-1} h(r) \le \delta_{E_r}(p^r) \le \beta h(r)$ \ \ \ \ $\beta^{-1} f(r) \le \delta_{E_r}(q^r) \le \beta f(r)$

\vspace{2mm}

$(v_b)$ The functions $h$ and $f$ satisfy the following Condition, called Condition(**): \vspace{2mm}

\vskip 0,1cm
\noindent Condition(**) {\sl ``$f=g^{-1}$ (resp. $h=k^{-1}$) where $g$ (resp. $k$) is a strictly increasing (resp. decreasing) convex function of class $\mathcal C^\infty$ defined on $[0,\varepsilon]$ (resp. $[-\epsilon,0]$) for some $\varepsilon > 0$ and such that $g^{(l)}(0) = 0$ (resp. $k^{(l)}(0)=0$) for every nonnegative integer $l$''.}

\vskip 0,2cm
We have:
\begin{lemma}\label{doubling-lem}
For every $0 < r < r_0$, for every $p \in D_r$ and for every  $v \in \mathbb C^n$ we have $K_{D_{2r}}(p,v) \leq 2 K_D(p,v)$ where $D_r:=D \cap \{z \in \mathbb C^n : Re(z_n) < r\}$.
\end{lemma}

\noindent{\bf Proof of Lemma~\ref{doubling-lem}.} Assume that $\delta_D(p,v) > \delta_{D_{2r}}(p,v)$. Since ${\partial D}_{2r}= ({\partial D} \cap \{z \in \mathbb C^n : Re(z_n) \leq 2r\}) \cup (\overline{D} \cap \{z \in \mathbb C^n : Re(z_n)=2r\}$  then there is $p \in D_r$, $v \in \mathbb C^n$ and $\lambda \in \mathbb C$ such that

$$
\delta_D(p,v) > dist_{eucl}(p,p+\lambda v)
$$

with $Re(p_n+\lambda v_n) = 2r$ (meaning $\delta_D(p,v) > \delta_{D_{2r}}(p,v)$).

Since $Re(p_n) < r$ then $Re(\lambda v_n) > r$ and so
$$
Re(p_n-\lambda v_n) < 0.
$$
In particular there is $0 < |\lambda '| < |\lambda|$ such that $Re(p_n+\lambda 'v_n) = 0$, meaning that the point $p+\lambda ' v \not\in D$. Hence $\delta_D(p,v) = \delta_{D_{2r}}(p,v)$ which is a contradiction. Hence

$$
\delta_D(p,v) = \delta_{D_{2r}}(p,v).
$$
It follows now from Proposition A that
$$
K_{D_{2r}}(p,v) \leq \frac{\|v\|}{\delta_{D_{2r}}(p,v)} = \frac{\|v\|}{\delta_D(p,v)} \leq 2 K_D(p,v).
$$
\qed

\vskip 0,2cm
For $0 < r < r_0$ we denote by:

\vskip 0,1cm
\noindent $\bullet$ $\gamma^r : [t_0,t_1] \rightarrow D$ a geodesic satisfying $\gamma^r(t_0) = p^r$, $\gamma^r(t_1) = q^r$,

\vskip 0,1cm
\noindent $\bullet$ $r':=\sup_{t \in [t_0,t_1]}Re(\gamma^r_n(t))$ where $\gamma^r:(\gamma^r_1,\dots,\gamma^r_n)$,

\vskip 0,1cm
\noindent $\bullet$ $\pi_1^{2r'}: \mathbb C^n \rightarrow E_{2r'}=\mathbb C \times \{(0,\dots,0,2r')\}$ denote the holomorphic projection onto the first factor,

\vskip 0,1cm
\noindent $\bullet$ $l_{p^r}^1$ (resp. $l_{q^r}^1$) the real line joining $p^r$ (resp. $q^r$) to $p^{2r'}=\pi_1^{2r'}(p^r)$ (resp. $q^{2r'}=\pi_1^{2r'}(q^r)$), contained in the real  2-plane $span_{\mathbb R}(e^1,e^n)$,

\vskip 0,1cm
\noindent $\bullet$ $\sigma^{2r'}$ the geodesic in $(E_{2r'},d^K_{E_{2r'}})$ joining $\pi_1^{2r'}(p^r)$ to $\pi_1^{2r'}(q^r)$.

\vskip 0,5cm

\centerline{\begin{picture}(0,0)%
\includegraphics{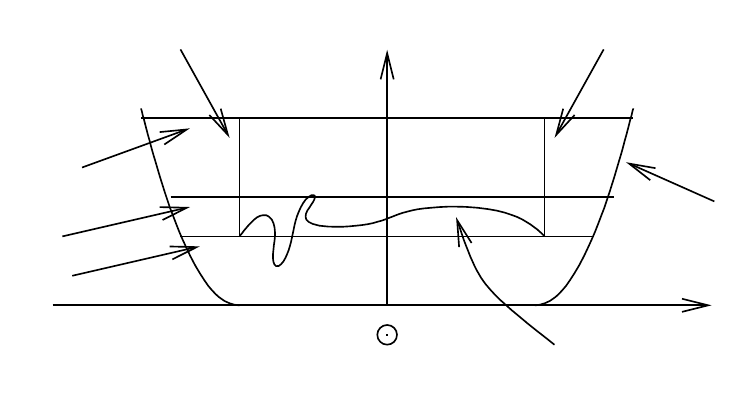}%
\end{picture}%
\setlength{\unitlength}{4144sp}%
\begingroup\makeatletter\ifx\SetFigFont\undefined%
\gdef\SetFigFont#1#2#3#4#5{%
  \reset@font\fontsize{#1}{#2pt}%
  \fontfamily{#3}\fontseries{#4}\fontshape{#5}%
  \selectfont}%
\fi\endgroup%
\begin{picture}(3360,1815)(706,-1840)
\put(1711,-1276){\makebox(0,0)[lb]{\smash{{\SetFigFont{12}{14.4}{\familydefault}{\mddefault}{\updefault}{\color[rgb]{0,0,0}$p^ r$}%
}}}}
\put(3061,-1276){\makebox(0,0)[lb]{\smash{{\SetFigFont{12}{14.4}{\familydefault}{\mddefault}{\updefault}{\color[rgb]{0,0,0}$q^ r$}%
}}}}
\put(2566,-286){\makebox(0,0)[lb]{\smash{{\SetFigFont{12}{14.4}{\familydefault}{\mddefault}{\updefault}{\color[rgb]{0,0,0}$Re(z_n)$}%
}}}}
\put(4051,-1411){\makebox(0,0)[lb]{\smash{{\SetFigFont{12}{14.4}{\familydefault}{\mddefault}{\updefault}{\color[rgb]{0,0,0}$z'$}%
}}}}
\put(3241,-1726){\makebox(0,0)[lb]{\smash{{\SetFigFont{12}{14.4}{\familydefault}{\mddefault}{\updefault}{\color[rgb]{0,0,0}$\gamma^ r$}%
}}}}
\put(2251,-1771){\makebox(0,0)[lb]{\smash{{\SetFigFont{12}{14.4}{\familydefault}{\mddefault}{\updefault}{\color[rgb]{0,0,0}$Im(z_n)$}%
}}}}
\put(3961,-1141){\makebox(0,0)[lb]{\smash{{\SetFigFont{12}{14.4}{\familydefault}{\mddefault}{\updefault}{\color[rgb]{0,0,0} $D_{2r'}$}%
}}}}
\put(721,-1366){\makebox(0,0)[lb]{\smash{{\SetFigFont{12}{14.4}{\familydefault}{\mddefault}{\updefault}{\color[rgb]{0,0,0}$E_r$}%
}}}}
\put(721,-1186){\makebox(0,0)[lb]{\smash{{\SetFigFont{12}{14.4}{\familydefault}{\mddefault}{\updefault}{\color[rgb]{0,0,0}$E_{r'}$}%
}}}}
\put(721,-871){\makebox(0,0)[lb]{\smash{{\SetFigFont{12}{14.4}{\familydefault}{\mddefault}{\updefault}{\color[rgb]{0,0,0}$E_{2r'}$}%
}}}}
\put(1351,-196){\makebox(0,0)[lb]{\smash{{\SetFigFont{12}{14.4}{\familydefault}{\mddefault}{\updefault}{\color[rgb]{0,0,0}$l^ 1_{p^ r}$}%
}}}}
\put(3511,-196){\makebox(0,0)[lb]{\smash{{\SetFigFont{12}{14.4}{\familydefault}{\mddefault}{\updefault}{\color[rgb]{0,0,0}$l^ 1_{q^ r}$}%
}}}}
\end{picture}%
}

\vskip 0,5cm

\noindent Then we have:
\begin{proposition}\label{quasi2-prop}
For $0 < r < r_0$ the real curve $\sigma^{r,r'}:=l_{q^r}^1 \cap \sigma^{2r'} \cap l_{p^r}^1$ is a quasi-geodesic connecting $p^r$ to $q^r$. Moreover, if $A_r$ and $B_r$ denote the corresponding constants given by Definition~\ref{quasi-def} then there is a positive constant $c$ such that $A_r > c$ and $B_r > c$.
\end{proposition}
Before proving Proposition~\ref{quasi2-prop} we need to check that $\pi_1^{2r'}(D \cap \{z \in \mathbb C^n : Re(z_n) \leq r'\}) \subset D$ for $0 < r < r_0$. By Condition ({\bf 4}) this will be true if $2r'<\sqrt{r_0}$. This is a consequence of the following lemma:

\begin{lemma}\label{geo-lem}
There exists $d > \frac{1}{2}$ such that for every $0 < r < r_0$ and for every $t \in [t_0,t_1]$ we have:
$$
Re(\gamma^r_n(t)) \leq r^d.
$$
Hence we have: $r' < \sqrt{r} < \sqrt{r_0}$.
\end{lemma}

\noindent{\bf Proof of Lemma~\ref{geo-lem}.} Assume to get a contradiction that there exists a sequence of points $t^\nu \in [t_0,t_1]$ and a sequence of numbers $d_\nu$ satisfying $\lim_{\nu \rightarrow \infty}d_\nu = \frac{1}{2}$ and $Re(\gamma^r_n(t^\nu)) \geq r^{d_\nu}$. Then according to the Boxing Lemma (Lemma~\ref{boxing-lem}) we have :
\begin{equation}\label{geo1-eq}
\forall \nu > 0,\ l_D^K(\gamma([t_0,t^\nu]) \geq \frac{1}{2}\left|\log\left(\frac{r^{d_\nu}}{r}\right)\right| = -\frac{1-d_\nu}{2}\log(r).
\end{equation}

Consider now the real line $L^r : [0,1] \rightarrow E_r \subset D$ given by $L^r(t)=(1-t)p^r+tq^r$. The Kobayashi length of $L^r$ can be estimated by :

\begin{align} \notag
 l^K_{E_r}(L^r) & \le   \int_0^1 \frac {\vert { (L^r)'(t)} \vert}{\delta_{E_r}(L^r(t))} dt \\ \notag
                &  \le  2{\alpha}^{-1} a \int_0^1 \frac {dt}{\delta_{E_r}(L^r(t), L^r)}  \\ \notag
                &  = 2{\alpha}^{-1} a \int_0^1 \frac {dt}{min \{a+f(r)-L_1^r(t), \ L_1^r(t)+a+h(r)\}}  \\ \notag
                &  \le 2{\alpha}^{-1} a \int_0^1 \frac {dt}{a+f(r)-L_1^r(t)} + 2{\alpha}^{-1} a \int_0^1 \frac {dt}{ L_1^r(t)+a+h(r)}  \\ \notag
\end{align}
where the second inequality is just Condition ({\bf 3}). The first integral above is:

$$ \int_0^1 \frac {dt}{a+f(r)-L_1^r(t)} = \int_0^1 \frac {dt}{2a+f(r)-2at} = -(2a)^{-1} log \Bigl ( \frac {f(r)}{2a+f(r)} \Bigr )$$
The second integral involving $h$ can be calculated in the same manner. Since $f$ and $h$ are increasing, we get
$$ l^K_{E_r}(L^r) \le -{\alpha}^{-1} log (f(r)) -{\alpha}^{-1} log (h(r)) +B$$
where $B= -{\alpha}^{-1}log (2a+f(r_0))-{\alpha}^{-1}log (2a+h(r_0))$.

Hence,if $0 < r < r_0,$
\begin{equation}\label{geo2-eq}
d^K_D(p^r,q^r) \leq d^K_{E_r}(p^r,q^r) \leq l^K_{E_r}(L^{r}) \leq C ( |log(f(r))| + | log (h(r))|) +B.
\end{equation}
where $C$ and $B$ do not depend on $r$.


Moreover we have:
\begin{lemma}\label{funct-lem}
Let $f$ and $h$ satisfy Condition (**). Then we have:
$$
\lim_{r \rightarrow 0}\left|\frac{\log f(r)}{\log r}\right| =\lim_{r \rightarrow 0}\left|\frac{\log h(r)}{\log r}\right| 0.
$$
\end{lemma}

\noindent{\bf Proof of Lemma~\ref{funct-lem}.} We prove this for $f$. Let $a(r)=r^l$. Then $a^{(l)}(0) > 1$ and so there exists $0 < \varepsilon_l < \varepsilon$ such that $g^{(l)}(r) < a^{(l)}(r)$ for every $0 < r < \varepsilon_l$. By integrating $l$-times that inequality we obtain:
$$
\forall 0 < r < \varepsilon_l,\ g(r) < r^l
$$
and so
$$
\forall 0 < r < \varepsilon_l,\ 1 > f(r) > r^{1/l}.
$$
Hence for every positive integer $l$, $|\log(f(r))| < \frac{1}{l}|\log r|$ for every $0 < r < \varepsilon_l$. This implies the desired statement. \qed

\vskip 0,2cm
Now the statement of Lemma~\ref{geo-lem} follows from (\ref{geo1-eq}), (\ref{geo2-eq}), Lemma \ref{funct-lem}
and the inequality $d^K_D(p^r,q^r) = l^K_D(\gamma^r) \ge l_D^K(\gamma^r([t_0,t^\nu])$. \qed

\vskip 0,2cm
\noindent{\bf Proof of Proposition~\ref{quasi2-prop}.}
To get the result we need to compare the Kobayashi lengths of the geodesic $\gamma^r$ joining $p^r$ to $q^r$ and of $\sigma^{r,r'}$.

We first observe that according to Lemma~\ref{doubling-lem} we have:
$$
l^K_D(\gamma^r) \geq \frac{1}{2}l^K_{D_{2r'}}(\gamma^r).
$$
Hence, since $\pi_1^{2r'}(D_{2r'}) \subset E_{2r'}$ by Condition 4 on $D$ and since $\sigma^{2r'}$ is a geodesic for $(E_{2r'},d^K_{E_{2r'}})$ it follows from the Schwarz Lemma (Decreasing property of the Kobayashi metric) that
$$
l^K_D(\sigma^{2r'}) \leq l^K_{E_{2r'}}(\sigma^{2r'}) \leq l^K_{E_{2r'}}(\pi_1^{2r'}(\gamma^r)) \leq l^K_{D_{2r'}}(\gamma^r) \leq 2l^K_D(\gamma^r).
$$
Moreover since by definition the curve $l_{p^r}^1$ is contained in the convex set $D^r_n:=D \cap (\{(p^r_1,0,\dots,0)\} \times \mathbb C)$ where $(p^r_1,0,\dots,0) \in \mathbb C^{n-1}$ then
$l^K_D(l^1_{p^r}) \leq l^K_{D^r_n}(l^1_{p^r})$.
However since by Condition 2 on $D$ we have $\delta_D(z) \geq \frac{Re(z_1)}{2}$ for every $z \in l¹_{p^r}$ then
$$
l^K_{D^r_n}(l^1_{p^r}) \leq 2 \left|\log \left(\frac{2r'}{r}\right)\right|.
$$
In a similar way we have:
$$
l^K_{D^r_n}(l^1_{q^r}) \leq 2 \left|\log \left(\frac{2r'}{r}\right)\right|.
$$


Finally it follows from Condition 1 on $D$ that there exist $\alpha, \beta, R>0$ such that $D \subset D_{R,\alpha,\beta}=\Delta_R^{n-1} \times R_{\alpha,\beta}$ where $R_{\alpha,\beta}:=\{\zeta \in \mathbb C : 0 < Re(\zeta) < \alpha,\ -\beta < Im(\zeta) < \beta\}$. Then according to the Boxing Lemma (Lemma~\ref{boxing-lem}) we have:
$$
l^K_D(\gamma^r) \geq C \left|\log\left(\frac{2r'}{r}\right)\right|+D
$$
where $C$ and $D$ are two positive constants independent of $r$ and $r'$.

Hence we proved that there is a positive constant $C'$ such that $l^K_D(\sigma^{r,r'}) \leq C' l^K_D(\gamma^r)$.
This gives the desired statement. \qed

\section{Estimates on quasi-triangles}
The aim of this section is to prove the following:

\begin{proposition}\label{nonhyp-prop}
Assume that $D$ satisfies all the properties of Section 5. Let $z^0:(z^0_1,\dots,z^0_n) \in D$ with $z^0_1 \in \mathbb R$ and let for $0 < r < r_0$, $p^r$ and $q^r$ be as in Section 5. Let $\sigma^{r,r'}$ be the quasi-geodesic joining $p^r$ to $q^r$ and defined in Proposition~\ref{quasi2-prop}. Let $l^{z^0,p^r}$ (resp. $l^{z^0,q^r}$) be the segments joining $z^0$ to $p^r$ (resp. to $q^r$).
Then there is a point $z^{r,r'} \in \sigma^{r,r'}$ such that
$$
\lim_{r \rightarrow 0}d^K_D(z^{r,r'},l^{z^0,p^r} \cup l^{z^0,q^r}) = +\infty.
$$
\end{proposition}

\noindent{\bf Proof of Proposition~\ref{nonhyp-prop}.} Since $D$ is bounded, smooth and convex, then there exists $\alpha_0>0$ such that every real line through $z^0$ intersects ${\partial D}$ with some angle $\alpha > \alpha_0$. Hence it follows from Proposition~\ref{quasi1-prop} that all such lines are quasi-geodesics with constants $A$ and $B$, given by Definition~\ref{quasi-def}, depending only on $z^0$ and $D$. It also follows from Proposition~\ref{quasi2-prop} that we choose $A$ and $B$ and $r_0>0$ such that $\sigma^{r,r'}$ is a $(A,B)$ quasi-geodesic for every $0 < r < r_0$ (we recall that $r' < \sqrt{r}$).

Moreover we have:
$$
{\partial D}_{\sqrt{2r'}}\cap l^{z^0,p^r}=\{p^{\sqrt{2r'}}=(p^{\sqrt{2r'}}_1,\dots,p^{\sqrt{2r'}}_n)\}
$$
with $p^{\sqrt{2r'}}_1 \in \mathbb R$
and
$$
{\partial D}_{\sqrt{2r'}}\cap l^{z^0,q^r}=\{q^{\sqrt{2r'}}=(q^{\sqrt{2r'}}_1,\dots,q^{\sqrt{2r'}}_n)\}
$$
with $q^{\sqrt{2r'}}_1 \in \mathbb R$.

Define the real hyperplanes
$$
H_{p^{\sqrt{2r'}}}:=\{z \in \mathbb C^n : Re(z_1) < p^{\sqrt{2r'}}_1\}
$$
and
$$
H_{q^{\sqrt{2r'}}}:=\{z \in \mathbb C^n : Re(z_1) > q^{\sqrt{2r'}}_1\}
$$

\vskip 0,5cm
\centerline{\begin{picture}(0,0)%
\includegraphics{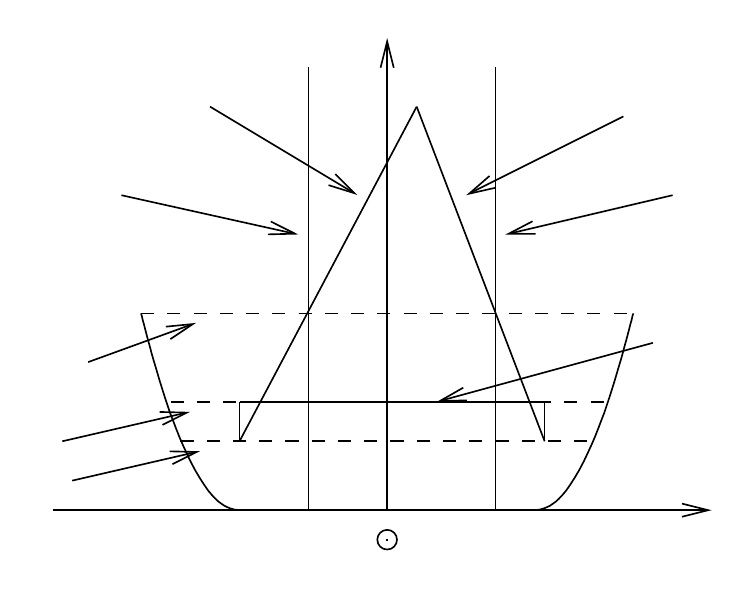}%
\end{picture}%
\setlength{\unitlength}{4144sp}%
\begingroup\makeatletter\ifx\SetFigFont\undefined%
\gdef\SetFigFont#1#2#3#4#5{%
  \reset@font\fontsize{#1}{#2pt}%
  \fontfamily{#3}\fontseries{#4}\fontshape{#5}%
  \selectfont}%
\fi\endgroup%
\begin{picture}(3360,2748)(706,-1840)
\put(1711,-1276){\makebox(0,0)[lb]{\smash{{\SetFigFont{12}{14.4}{\familydefault}{\mddefault}{\updefault}{\color[rgb]{0,0,0}$p^ r$}%
}}}}
\put(3061,-1276){\makebox(0,0)[lb]{\smash{{\SetFigFont{12}{14.4}{\familydefault}{\mddefault}{\updefault}{\color[rgb]{0,0,0}$q^ r$}%
}}}}
\put(4051,-1411){\makebox(0,0)[lb]{\smash{{\SetFigFont{12}{14.4}{\familydefault}{\mddefault}{\updefault}{\color[rgb]{0,0,0}$z'$}%
}}}}
\put(2251,-1771){\makebox(0,0)[lb]{\smash{{\SetFigFont{12}{14.4}{\familydefault}{\mddefault}{\updefault}{\color[rgb]{0,0,0}$Im(z_n)$}%
}}}}
\put(721,-1366){\makebox(0,0)[lb]{\smash{{\SetFigFont{12}{14.4}{\familydefault}{\mddefault}{\updefault}{\color[rgb]{0,0,0}$E_r$}%
}}}}
\put(721,-1186){\makebox(0,0)[lb]{\smash{{\SetFigFont{12}{14.4}{\familydefault}{\mddefault}{\updefault}{\color[rgb]{0,0,0}$E_{2r'}$}%
}}}}
\put(3736,-691){\makebox(0,0)[lb]{\smash{{\SetFigFont{12}{14.4}{\familydefault}{\mddefault}{\updefault}{\color[rgb]{0,0,0}$\sigma^ {r,r'}$}%
}}}}
\put(3466,434){\makebox(0,0)[lb]{\smash{{\SetFigFont{12}{14.4}{\familydefault}{\mddefault}{\updefault}{\color[rgb]{0,0,0}$l^ {z^ 0,q^ r}$}%
}}}}
\put(1396,479){\makebox(0,0)[lb]{\smash{{\SetFigFont{12}{14.4}{\familydefault}{\mddefault}{\updefault}{\color[rgb]{0,0,0}$l^ {z^ 0,p^ r}$}%
}}}}
\put(991,119){\makebox(0,0)[lb]{\smash{{\SetFigFont{12}{14.4}{\familydefault}{\mddefault}{\updefault}{\color[rgb]{0,0,0}$H_{p^ {\sqrt{2r'}}}$}%
}}}}
\put(3736, 74){\makebox(0,0)[lb]{\smash{{\SetFigFont{12}{14.4}{\familydefault}{\mddefault}{\updefault}{\color[rgb]{0,0,0}$H_{q^ {\sqrt{2r'}}}$}%
}}}}
\put(2521,-421){\makebox(0,0)[lb]{\smash{{\SetFigFont{12}{14.4}{\familydefault}{\mddefault}{\updefault}{\color[rgb]{0,0,0}$z^ {r,r'}$}%
}}}}
\put(2476,749){\makebox(0,0)[lb]{\smash{{\SetFigFont{12}{14.4}{\familydefault}{\mddefault}{\updefault}{\color[rgb]{0,0,0}$Re(z_n)$}%
}}}}
\put(2656,479){\makebox(0,0)[lb]{\smash{{\SetFigFont{12}{14.4}{\familydefault}{\mddefault}{\updefault}{\color[rgb]{0,0,0}$z^ 0$}%
}}}}
\put(721,-826){\makebox(0,0)[lb]{\smash{{\SetFigFont{12}{14.4}{\familydefault}{\mddefault}{\updefault}{\color[rgb]{0,0,0}$E_{\sqrt{2r'}}$}%
}}}}
\end{picture}%
}
\vskip 0,5cm

We first point out that by construction every point $z\in \sigma^{r,r'}$ satisfies $Re(z_1) \leq -2r'$. Hence it follows from Lemma
~\ref{boxing-lem} that
$$
d^K_D(z,l^{z^0,p^r} \cap D_{\sqrt{2r'}}^c) \geq \frac{1}{2}\left|\log\left(\frac{\sqrt{2r'}}{2r'}\right)\right|.
$$
Since $r'<\sqrt{r}$ according to Lemma~\ref{geo-lem}, we have:
\begin{equation}\label{last-eq}
\lim_{r \rightarrow 0}\inf_{z \in \sigma^{r,r'}}d^K_D(z,l^{z^0,p^r}\cap D_{\sqrt{2r'}}^c) = +\infty
\end{equation}
where $D_{\sqrt{2r'}}^c$ denotes $D \backslash D_{\sqrt{2r'}}$.

Let $z^{r,r'}$ be any point on $\sigma^{r,r'} \cap E_{2r'}$ satisfying $Re(z^{r,r'}_1)=0$. Such a point exists since by construction $\sigma^{r,r'}=l^{p^r}_1 \cup \sigma^{2r'} \cup l^{q^r}_1$ and $\sigma^{2r'}$ joins the two points $\pi_1^{2r'}(p^r)$ and $\pi_1^{2r'}(q^r)$ with $p^r_1<0$ and $q^r_1 > 0$.

According to (\ref{last-eq}), to prove Proposition~\ref{nonhyp-prop} we must prove that $\lim_{r \rightarrow 0}d^K_D(z^{r,r'},l^{z^0,p^r} \cap D_{\sqrt{2r'}})=+\infty$.
Since every point $x \in l^{z^0,p^r} \cap D_{\sqrt{2r'}}$ satisfies $Re(x_1)<p^{\sqrt{2r'}}_1$ we just need the following condition:

\begin{equation}\label{p-eq}
\lim_{r \rightarrow 0}d^K_D(z^{r,r'},H_{p^{\sqrt{2r'}}})=+\infty.
\end{equation}

It follows from Lemma~\ref{geo-lem} that any point $z$ on a geodesic $\tilde \gamma^{r,r'}$ in $(D,d^K_D)$ joining $z^{r,r'}$ to $H_{p^{\sqrt{2r'}}}$ will satisfy $Re(z_n) > -(2r')^{1/4}$.
Hence according to Lemma~\ref{doubling-lem} we have:
$$
l^K_D(\tilde \gamma^{r,r'}) \geq \frac{1}{2}l^K_{D_{2(2r')^{1/4}}}(\tilde \gamma^{r,r'}) \geq \frac{1}{2}l^K_{E_{2(2r')^{1/4}}}(\pi^{2(2r')^{1/4}}(\tilde \gamma^{r,r'}),
$$
the last inequality coming from the Schwarz Lemma.

For convenience we set $s:=2(2r')^{1/4}$ and we assume that $\tilde \gamma^{r,r'} : [0,1] \rightarrow D$. Finally we set $\gamma:=\pi_1^s(\tilde \gamma^{r,r'})$. Then $Re(\gamma_1(0))=0$ and $Re(\gamma_1(1)))=p^{\sqrt{2r'}}_1$ where $\gamma=(\gamma_1,\dots,\gamma_n)$.

Finally $\delta_{E_s}(\gamma(t)) \leq 1-Re(\gamma_1(t))$. So
$$
\begin{array}{lllll}
l^K_{E_s}(\gamma) & \geq & \displaystyle \frac{1}{2}\int_0^1\frac{|Re(\gamma_1'(t))|}{1-Re(\gamma_1(t))}dt & \geq & \displaystyle \frac{1}{2} \left|\int_0^1\frac{Re(\gamma_1'(t))}{1-Re(\gamma_1(t))}dt\right|\\
& & & & \\
& & & \geq & \displaystyle \frac{-\log(1-p^{\sqrt{2r'}}_1)}{2} \rightarrow_{r \rightarrow 0}+\infty.
\end{array}
$$
This gives:
$$
\lim_{r \rightarrow 0}l^K_D(\tilde \gamma^{r,r'})=+\infty
$$
or equivalently
$$
\lim_{r \rightarrow 0}d^K_D(z^{r,r'},H_{p^{\sqrt{2r'}}} \cap D_{\sqrt{2r'}}^c)=+\infty.
$$
This proves Condition~(\ref{p-eq}).
In a similar way we obtain:

$$
\lim_{r \rightarrow 0}d^K_D(z^{r,r'},H_{q^{\sqrt{2r'}}}\cap D_{\sqrt{2r'}}^c)=+\infty.
$$

This gives the desired result. \qed




\section{Gromov-hyperbolicity of complex ellipsoids}
Consider the following {\it complex ellipsoid in} $\C^n$:
$$D_p = \{(z',z_n) \in \C^n \ \vert \ \vert z' \vert^2 + \vert z_n \vert^{2p} < \ 1\},$$
where $p \ge 1$. \vspace{2mm}

\begin{theorem}\label{CompEllip-thm2}
For $p \geq 1$ the complex ellipsoid $D_p$ is Gromov hyperbolic for the Kobayashi distance.
\end{theorem}

The proof of Theorem~\ref{CompEllip-thm2} will be a consequence of the following facts.

\begin{lemma}\label{gh}
 Let $D \subset \C^n$ be a bounded domain.

(1) Let $g$ be a K\"ahler metric on $D$ with uniformly negative holomorphic
sectional curvature. Then there exists $C_1 >0$ such that
$$  \sqrt {g(v,v)} \ \le \ C_1 K_D(p,v) $$
for all $p \in D$, $v \in \C^n$.

(2) Let $g$ be a complete K\"ahler metric on $D$ whose Ricci curvature is uniformly bounded from below.
Then there exists $C_2 >0$ such that
$$   C_2 C_D(p,v) \ \le \ \sqrt  {g(v,v)} $$
for all $p \in D$, $v \in \C^n$.
\end{lemma}

{\bf Proof:} These are immediate consequences of the Yau-Schwarz Lemma ~\cite{yau}: Let $(M,g)$ and $(N,h)$ be K\"ahler manifolds. Assume that $(M,g)$ is complete and that there exists positive constants $A, B$ such that $Ric_g \ge -A$ and
$H_h \le -B$, where $Ric_g$ and $H_h$ denote the Ricci curvature and holomorphic sectional curvatures of $g$ and
$h$ respectively.  If $f:M \rt N$ is any holomorphic map then
$$f^\ast h \le \frac{A}{B} g.$$
To see $(1)$ let $p \in D, v \in \C^n, w \in \C$ and $f: \Delta \rt D$ with $f(0)=p$ and $df_0(w)=v$. By the Yau-Schwarz
Lemma (here we consider $\Delta$ endowed with the Poincare metric and note that the Poincar\'e metric on $\Delta$ coincides with the Euclidean metric at the origin), we have
$$ \sqrt  {g(df_0(w), df_0(w))}  \le  C_1 \vert w \vert $$
where $-C_1<0$ is a lower bound on the holomorphic sectional curvatures of $(M,g)$. Since this inequality
holds for all $w$ satisfying $df_0(w)=v$ we get $(1)$.

For $(2)$, consider a holomorphic map $f: D \rt \Delta$ with $f(p)=0$ and let $df_p(v)=w$ where $v \in \C^n$. Again,
endowing $\Delta$ with the Poincare metric and applying the Yau-Schwarz Lemma we get
$$ \vert df_p(v) \vert \le  A \sqrt  {g(v, v)} $$
where $A$ is a lower bound for the Ricci curvature of $(D,g)$. As earlier this gives $(2)$ with $C_2= A^{-1}$. \hfill $\square$
\vspace{2mm}

{\bf Remark:} It follows from the work of K. Azukawa and M. Suzuki ~\cite{azu} that the  {\it holomorphic} sectional curvatures of the Bergman metric on $D$ are bounded between two negative constants. Hence, by Lemma \ref{gh}, it follows that the Bergman metric is bi-Lipschitz to the Kobayashi metric as well.

Recall that two Finsler metrics $F_1$ and $F_2$ on a domain $D$ are said to be bi-Lipschitz if there exists a positive
constant $A$ such that
$$ A^{-1}F_1(p,v) \ \le \ F_2(p,v) \ \le \ AF_1(p,v)$$
for all $(p,v) \in D \times \C^n$.  Similarly two distance functions $d_1$ and $d_2$ on $D$ are bi-Lipschitz if there exists a positive constant $A$ such that
$$ A^{-1}d_1(p,q) \ \le \ d_2(p,q) \ \le \ Ad_1(p,q)$$
for all $p, q$ in $D$. Of course, if the distance functions arise from Finsler metrics
then the bi-Lipschitz equivalence of the metrics implies that of the distance functions.  \vspace{2mm}

\begin{lemma}\label{jk}
The Kobayashi and the K\"ahler-Einstein metrics are bi-Lipschitz on $D$.
\end{lemma}
{\bf Proof:}  Since $D$ is convex, it follows by the work of L. Lempert ~\cite{lem1}  that $C_D=K_D$. The result of Bland that the sectional curvatures of the K\"ahler-Einstein metric are bounded between two negative constants and
Lemma \ref{gh} then complete the proof. \hfill $\square$ \vspace{2mm}

\noindent{\bf Proof of Theorem~\ref{CompEllip-thm2}.} One knows, by the Toponogov Comparison Theorem in Riemannian geometry \cite{ch-eb}, that a simply-connected complete Riemannian manifold with sectional curvature bounded from above by a negative constant is Gromov-hyperbolic. Also, by computations due to
J. Bland ~\cite{bla}, the  K\"ahler-Einstein metric on $D$ has sectional curvatures bounded between two negative constants.
Since $D$ is simply-connected then $D$ is Gromov-hyperbolic for the distance function arising from the K\"ahler-Einstein metric. By Lemma \ref{jk} this distance is bi-Lipschitz to the Kobayashi distance. Since Gromov-hyperbolicity is
preserved by the notion of bi-Lipschitz equivalence (in fact, it is preserved under the weaker condition
of {\it quasi-isometric} equivalence) the theorem follows. \hfill $\square$ 

\vskip 0,5cm
{\small
\noindent Herv\'e Gaussier\\
(1) UJF-Grenoble 1, Institut Fourier, Grenoble, F-38402, France\\
(2) CNRS UMR5582, Institut Fourier, Grenoble, F-38041, France\\
{\sl E-mail address} : herve.gaussier@ujf-grenoble.fr\\
\\
Harish Seshadri\\
Department of Mathematics, Indian Institute of Science, Bangalore 560012, India\\
{\sl E-mail address} : harish@math.iisc.ernet.in
}

\end{document}